\documentclass[12pt,reqno]{amsart}
\usepackage{fullpage}
\usepackage{amsfonts}
\usepackage{times} 
\usepackage{graphicx}

\newcommand{\ep}{\varepsilon}
\newcommand{\M}{\mathcal{M}}
\newcommand{\I}{\mathcal{I}}
\newcommand{\K}{k}
\newcommand{\Q}{{\mathbb Q}}
\newcommand{\R}{{\mathbb R}}
\newcommand{\Z}{{\mathbb Z}}
\renewcommand{\v}{\mathbf{v}}
\newcommand{\va}{{\mathbf a}}
\newcommand{\ve}{{\mathbf e}}
\newcommand{\vm}{{\mathbf m}}

\newtheorem{lemma}{Lemma}

\newtheorem*{theorem}{Theorem}
{
\theoremstyle{remark}

\newtheorem*{remark}{Remark}
}

\begin{document}

\title{Almost all integer matrices have no integer eigenvalues}
\author{Greg Martin and Erick B. Wong}
\address{Department of Mathematics \\ University of British Columbia \\ Room
121, 1984 Mathematics Road \\ Canada V6T 1Z2}
\email{gerg@math.ubc.ca and erick@math.ubc.ca}
\maketitle

\section{Introduction}
In a recent issue of this {\sc Monthly}, Hetzel, Liew, and Morrison \cite{HLM} pose a rather natural question: {\em what is the probability that a random $n \times n$ integer matrix is diagonalizable over the field of rational numbers?}

Since there is no uniform probability distribution on $\Z$, we need to exercise some care in interpreting this question.  Specifically, for an integer $\K\ge 1$, let $\I_\K = \{ -\K, -\K+1, \dots, \K-1, \K \}$ be the set of integers with absolute value at most~$\K$.  Since $\I_\K$ is finite, we are free to choose each entry of an $n \times n$ matrix $M$ independently and uniformly at random from $\I_\K$, with each value having probability $1/(2\K+1)$.  The probability that $M$ has a given property, such as being diagonalizable over $\Q$, is then a function of the parameter $\K$; we consider how this function behaves as $\K \to \infty$.  In particular, if the probability converges to some limiting value then it is natural to think of this limit as the ``probability'' that a random integer matrix has that property.

We refer the reader to the article of Hetzel {\it et al} for an interesting discussion of some of the issues raised by this interpretation of probability as a limit of finitary probabilities (in doing so we lose countable additivity and hence the  measure-theoretic foundation of modern probability theory after Kolmogorov).  From a pragmatic viewpoint, this cost is outweighed by the fact that many beautiful number-theoretic results are most naturally phrased in the language of probability: for instance, the celebrated Erd\H{o}s-Kac theorem \cite{EK} states that the number of prime factors of a positive integer $n$ behaves (in an appropriate limiting sense) as a normally-distributed random variable with mean and variance both equal to $\log \log n$. (In this article we always mean the natural logarithm when we write $\log$.)

For any given integers $n \ge 2$ and $\K \ge 1$, the set of random $n \times n$ matrices with entries in $\I_\K$ is a finite probability space; it will be convenient to compute probabilities simply by counting matrices, so we introduce some notation for them.  Let $\M_n(\K)$ denote the set of all $n\times n$ matrices whose entries are all in $\I_\K$; then we are choosing matrices uniformly from $\M_n(\K)$, which has cardinality exactly $(2\K+1)^{n^2}$.  The probability that a random matrix in $\M_n(\K)$ satisfies a particular property is simply the number of matrices in $\M_n(\K)$ with that property divided by $(2\K+1)^{n^2}$.

For a given integer $\lambda$, let $\M^\lambda_n(\K)$ denote the set of all matrices in $\M_n(\K)$ that have $\lambda$ as an eigenvalue.  Note that in particular, $\M^0_n(\K)$ is the subset of singular matrices in $\M_n(\K)$.  Likewise, we denote set of the matrices in $\M_n(\K)$ having at least one integer eigenvalue by $\M^\Z_n(\K) = \bigcup_{\lambda \in \Z} \M^\lambda_n(\K)$.  The probability that a random matrix in $\M_n(\K)$ has an integer eigenvalue is thus $|\M^\Z_n(\K)|/(2\K+1)^{n^2}$.

Our main result affirms and strengthens a conjecture made in~\cite{HLM}: {\em for any $n\ge 2$, the probability that a random $n \times n$ integer matrix has even a single integer eigenvalue is~0}.  We furthermore give a quantitative upper bound on the decay rate of the probability as $\K$ increases. It will be extremely convenient to use ``Vinogradov's notation'' to express this decay rate: we write $f(\K) \ll g(\K)$ if there exists a constant $C>0$ such that $|f(\K)| \le Cg(\K)$ for all values of $\K$ under consideration. Notice for example that if $f_1(\K) \ll g(\K)$ and $f_2(\K) \ll g(\K)$, then $f_1(\K) + f_2(\K) \ll g(\K)$ as well. If this constant can depend upon some auxiliary parameter such as $\ep$, then we write $f(\K) \ll_\ep g(\K)$; for example, for $\K\ge1$ it is true that $\log \K \ll_\ep \K^\ep$ for every $\ep>0$.

\begin{theorem}
Given any integer $n\ge2$ and any real number $\ep>0$, the probability that a randomly chosen matrix in $\M_n(\K)$ has an integer eigenvalue is $\ll_{n,\ep} 1/\K^{1-\ep}$.  In particular, the probability that a randomly chosen matrix in $\M_n(\K)$ is diagonalizable over the rational numbers is $\ll_{n,\ep} 1/\K^{1-\ep}$.
\label{theorem}
\end{theorem}

Given an integer matrix $M \in \M_n(\K)$, a necessary condition for it to be diagonalizable over $\Q$ is that all of its eigenvalues are rational.  Moreover, since the characteristic polynomial $\det(\lambda I - M)$ is monic with integer coefficients, the familiar ``rational roots theorem'' \cite[\S4.3]{N} implies that every rational eigenvalue of $M$ must be an integer. Hence any matrix that is diagonalizable over the rationals must certainly belong to $\M^\Z_n(\K)$, and so the second assertion of the theorem follows immediately from the first.

The special case $n=2$ of the theorem was obtained in \cite{HLM} and also earlier by Kowalsky \cite{Kow}.  Unravelling the $\ll$-notation, the theorem states that there exists a constant $C$, possibly depending on $n$ and $\ep$, such that $|\M^\Z_n(\K)|/|\M_n(\K)| \le C/\K^{1-\ep}$ for all $\K \ge 1$. Note that $|\M_n(\K)| \ll_n \K^{n^2}$ (with the implied constant highly dependent on $n$), and so the theorem also gives an upper bound for the number of matrices in $\M_n(\K)$ with at least one integer eigenvalue, namely
\[
|\M^\Z_n(\K)| \ll_{n,\ep} \K^{n^2-1+\ep}.
\]

The key tool used to establish the theorem is the following related estimate for the number of singular matrices in $\M_n(\K)$:

\begin{lemma}
Given any integer $n\ge2$ and any real number $\ep>0$, the probability that a random matrix in $\M_n(\K)$ is singular is $\ll_{n,\ep} 1/\K^{2-\ep}$. In other words, $|\M^0_n(\K)| \ll_{n,\ep} \K^{n^2-2+\ep}$.
\label{any.dimension}
\end{lemma}

This upper bound is essentially best possible when $n=2$: we can show that the cardinality of $\M^0_2(\K)$ is asymptotic to $(96/\pi^2)\K^2 \log \K$, so that the probability that a matrix in $\M_2(\K)$ is singular is asymptotic to $(6\log\K)/\pi^2\K^2$.  (We discuss the $2\times 2$ case in more detail in Section~\ref{2by2}.) However, we expect that the upper bound of Lemma~\ref{any.dimension} is not sharp for any $n>2$. The probability that a matrix has a row consisting of all zeros, or that two of its rows are identical, is $1/\K^n$ (up to a constant depending only on $n$), so the upper bound cannot be any smaller than this. It seems reasonable to conjecture that the true proportion of singular matrices in $\M_n(\K)$ decays as $1/\K^{n-\ep}$, but for higher dimensions we do not know what the correct analogy should be of the precise rate $(\log \K)/\K^2$ that holds in dimension~2.

We remark that in this paper, we are considering the behaviour of $|\M^0_n(\K)|/|\M_n(\K)|$ when $n$ is fixed and $\K$ increases.  It is perfectly natural to ask how this same probability behaves when $\K$ is fixed and $n$ increases, and there is a body of deep work on this problem.  Koml\'{o}s \cite{Kom} proved in 1968 that this probability converges to 0 as $n \to \infty$, answering a question of Erd\H{o}s.  In fact Koml\'{o}s's result holds, not just for the uniform distribution on $\I_\K$, but for an {\em arbitrary} non-degenerate distribution on $\R$. (A degenerate distribution is constant almost surely, so the non-degeneracy condition is clearly necessary.)  Slinko \cite{S} later established a quantitative decay rate of $\ll_\K 1/\sqrt{n}$.  An exponential decay rate of $(0.999)^n$ for the case of $\{\pm 1\}$-matrices was established by Kahn, Koml\'{o}s, and Szemer\'{e}di \cite{KKS} and improved to $(3/4)^n$ by Tao and Vu \cite{TV}.  Very recently Rudelson and Vershynin \cite{RV} have established exponential decay for a very wide class of distributions, including the uniform distribution on $\I_\K$.

\section{Determinants, singular matrices, and integer eigenvalues}

We begin by proving a lemma that we will use repeatedly in the proof of Lemma~\ref{any.dimension}.  It shows that the probability that a $2\times 2$ matrix is singular remains small (that is, the probability is $\ll \K^{-2+\ep}$ just as in Lemma~\ref{any.dimension}, albeit with a different implied constant) even if we choose the entries randomly from arbitrary arithmetic progressions of the same length as $\I_\K$.

\begin{lemma}
Fix positive real numbers $\alpha$ and $\ep$. Let $k$ be a positive integer, and let $L_1(x)$, $L_2(x)$, $L_3(x)$, and $L_4(x)$ be non-constant linear polynomials whose coefficients are integers that are $\ll_\alpha \K^\alpha$ in absolute value. Then the number of solutions to the equation
\begin{equation}
L_1(a)L_2(b) = L_3(c)L_4(d),
\label{linear.products}
\end{equation}
with all of $a$, $b$, $c$, and $d$ in $\I_\K$, is $\ll_{\alpha,\ep} \K^{2+\ep}$.
\label{not.many.solutions}
\end{lemma}

\begin{proof}
First we consider the solutions for which both sides of equation \eqref{linear.products} equal~0. In this case, at least two of the linear factors $L_1(a)$, $L_2(b)$, $L_3(c)$, and $L_4(d)$ equal~0. If, for example, $L_1(a)=0$ and $L_3(c)=0$ (the other cases are exactly the same), this completely determines the values of $a$ and $c$; since there are $2\K+1$ choices for each of $b$ and $d$, the total number of solutions for which both sides of equation \eqref{linear.products} equal~0 is $\ll \K^2$.

Otherwise, fix any values for $c$ and $d$ for which the right-hand side of equation \eqref{linear.products} is nonzero, a total of at most $(2\K+1)^2 \ll \K^2$ choices. Then the right-hand side is some nonzero integer that is at most $(\K\cdot \K^\alpha+\K^\alpha)^2 \le 4\K^{2+2\alpha}$ in absolute value, and $L_1(a)$ must be a divisor of that integer.

It is a well-known lemma in analytic number theory (see for instance \cite[p.~56]{Magic}) that for any $\delta>0$, the number of divisors of a nonzero integer $n$ is $\ll_\delta |n|^\delta$.
In particular, choosing $\delta = \ep/(2+2\alpha)$, the right-hand side of equation \eqref{linear.products} has $\ll_{\alpha,\ep} (4\K^{2+2\alpha})^{\ep/(2+2\alpha)} \ll_{\alpha,\ep} \K^\ep$ divisors to serve as candidates for $L_1(a)$; each of these completely determines a possibility for $a$ (which might not even be an integer). Then the possible values for $L_2(b)$ and hence $b$ are determined as well. We conclude that there are a total of $\ll_{\alpha,\ep} \K^{2+\ep}$ solutions to equation \eqref{linear.products} as claimed.
\end{proof}

\begin{remark}
It is not important that the $L_i$ be {\em linear} polynomials: the above proof works essentially without change for any four non-constant polynomials of bounded degree. We will not need such a generalization, however, as the determinant of a matrix depends only linearly on each matrix element.
\end{remark}

The next ingredient is a curious determinantal identity, which was classically known but at present appears to have fallen out of common knowledge.  Before we can state this identity, we need to define some preliminary notation.  For the remainder of this section, capital letters will denote matrices, boldface lowercase letters will denote column vectors, and regular lowercase letters will denote scalars.

Let $I_n$ denote the $n\times n$ identity matrix, and let $\ve_j$ denote the $j$th standard basis vector (that is, the $j$th column of $I_n$). Let $M$ be an $n\times n$ matrix, and let $\vm_j$ denote its $j$th column and $m_{ij}$ its $ij$th entry. Note that $M\ve_j = \vm_j$ by the definition of matrix multiplication.

Let $a_{ij}$ denote the $ij$th cofactor of $M$, that is, the determinant of the $(n-1)\times(n-1)$ matrix obtained from $M$ by deleting its $i$th row and $j$th column. Let $A = \mathop{\rm Adj}(M)$ denote the adjugate matrix of $M$, that is, the matrix whose $ij$th entry is $(-1)^{i+j}a_{ji}$.  It is a standard consequence of Laplace's determinant expansion \cite[\S4.III]{H} that $MA = (\det M)I_n$.  Finally, let $\va_j$ denote the $j$th column of $A$. Note that $M\va_j = (\det M)\ve_j$, since both sides are the $j$th column of $(\det M)I_n$.

\begin{lemma}
Fix an integer $n\ge3$. Given an $n\times n$ matrix $M$, let $a_{ij}$ denote the $ij$th cofactor of $M$. Also let $Z$ denote the $(n-2)\times(n-2)$ matrix obtained from $M$ by deleting the first two rows and first two columns, so that
\begin{equation}
M = \left( \begin{tabular}{cc|c} $m_{11}$ & $m_{12}$ & * \\ $m_{21}$ & $m_{22}$ & * \\ \hline $*$ & $*$ & $Z$ \end{tabular} \right).
\label{Z.form}
\end{equation}
Then $a_{11}a_{22} - a_{12}a_{21} = (\det M)(\det Z)$.
\label{determinant.identity}
\end{lemma}

It is important to note that when $\det Z\ne 0$, the cofactor $a_{11}$ is a linear polynomial in the variable $m_{22}$ with leading coefficient $\det Z$, while the cofactor $a_{22}$ is a linear polynomial in $m_{11}$ with leading coefficient $\det Z$ (and similarly for the pair $a_{12}$ and $a_{21}$).
For example, when $n=3$ the determinant of the $1\times1$ matrix $Z$ is simply the lower-right entry $m_{33}$ of $M$; the identity in question is thus
\begin{multline}
(m_{11}m_{33}-m_{13}m_{31})(m_{22}m_{33}-m_{23}m_{32})  \\
- (m_{12}m_{33}-m_{13}m_{32})(m_{21}m_{33}-m_{23}m_{31})
= m_{33} \det M.
\label{in.dimension.3}
\end{multline}
For any given dimension $n$, the assertion of Lemma \ref{determinant.identity} is simply some polynomial identity that can be checked directly; however, a proof that works for all $n$ at once requires a bit of cunning.

\begin{proof}
Define a matrix
\[
B = \big( \va_1 \ \va_2 \ \ve_3 \ \cdots \ \ve_n \big) = \left( \begin{tabular}{cc|c} $a_{11}$ & $-a_{21}$ & 0 \\ $-a_{12}$ & $a_{22}$ & 0 \\ \hline $*$ & $*$ & $I_{n-2}$ \end{tabular} \right),
\]
where $*$ represents irrelevant entries. Since $B$ is in lower-triangular block form, its determinant
\[
\det B = \det \left( \begin{matrix} a_{11} & -a_{21} \\ -a_{12} & a_{22} \end{matrix} \right) \cdot \det I_{n-2} = (a_{11}a_{22} - a_{12}a_{21})
\]
is easy to evaluate. Moreover,
\begin{align*}
MB &= \big( M\va_1 \ M\va_2 \ M\ve_3 \ \cdots \ M\ve_n \big) \\
&= \big( (\det M)\ve_1 \ (\det M)\ve_2 \ \vm_3 \ \cdots \ \vm_n \big) = \left( \begin{tabular}{cc|c} $\det M$ & 0 & $*$ \\ 0 & $\det M$ & $*$ \\ \hline 0 & 0 & Z \end{tabular} \right).
\end{align*}
Since $MB$ is in upper-triangular block form, its determinant $\det(MB) = (\det M)^2(\det Z)$ is also easy to evaluate. Using the identity $\det M\cdot \det B = \det(MB)$, we conclude
\[
(\det M)(a_{11}a_{22} - a_{12}a_{21}) = (\det M)^2(\det Z).
\]
Both sides of this last identity are polynomial functions of the $n^2$ variables $m_{ij}$ representing the entries of $M$. The factor $\det M$ on both sides is a nonzero polynomial, and hence it can be canceled to obtain $(\det M)(\det Z) = a_{11}a_{22} - a_{12}a_{21}$ as desired.
\end{proof}

\begin{remark}
This proof generalises readily to a similar statement for larger minors of the adjugate matrix $A$.
Muir's classic treatise on determinants \cite[Ch.~VI, \S175]{MM} includes this generalization of Lemma~\ref{determinant.identity} in a chapter wholly devoted to compound determinants (that is, determinants of matrices whose elements are themselves determinants).  The same result can also be found in Scott's reference of equally old vintage \cite[p.~62]{Sc}, which has been made freely available online by the Cornell University Library Historical Math collection.
\end{remark}

We are now ready to prove Lemma~\ref{any.dimension}. We proceed by induction on $n$, establishing both $n=2$ and $n=3$ as base cases.

{\bf Base case $n=2$}: The determinant of
$\big( \genfrac{}{}{0pt}{1}{m_{11}}{m_{21}} \genfrac{}{}{0pt}{1}{m_{12}}{m_{22}} \big)$
equals 0 precisely when $m_{11}m_{22} = m_{12}m_{21}$. By Lemma \ref{not.many.solutions}, there are $\ll_\ep \K^{2+\ep}$ solutions to this equation with the variables $m_{11}$, $m_{12}$, $m_{21}$, and $m_{22}$ all in $\I_\K$. This immediately shows that the number of matrices in $\M^0_2(\K)$ is $\ll_\ep \K^{2+\ep}$ as claimed.  Since $1/|\M_2(\K)| \ll \K^{-4}$, we see that the probability of a randomly chosen matrix from $\M_2(\K)$ being singular is $|\M^0_2(\K)|/|\M_2(\K)| \ll_\ep \K^{-2+\ep}$.

{\bf Base case $n=3$}: We first estimate the number of matrices in $\M^0_3(\K)$ whose lower right-hand entry $m_{33}$ is nonzero. Fix the five entries in the last row and last column of $M$, with $m_{33}\ne0$; there are a total of $2\K(2\K+1)^4 \ll \K^5$ possibilities. Using the identity \eqref{in.dimension.3}, we see that if $\det M=0$ then we must have
\[
(m_{11}m_{33}-m_{13}m_{31})(m_{22}m_{33}-m_{23}m_{32}) = (m_{12}m_{33}-m_{13}m_{32})(m_{21}m_{33}-m_{23}m_{31}).
\]
This equation is of the form $L_1(m_{11})L_2(m_{22}) = L_3(m_{12})L_4(m_{21})$, where the $L_i$ are non-constant linear polynomials whose coefficients are at most $\K^2$ in absolute value. (Note that we have used the fact that $m_{33}\ne0$ in asserting that the $L_i$ are non-constant.) Applying Lemma \ref{not.many.solutions} with $\alpha=2$, we see that there are $\ll_\ep \K^{2+\ep}$ solutions to this equation with $m_{11}$, $m_{12}$, $m_{21}$, and $m_{22}$ all in $\I_\K$. This shows that there are $\ll_\ep \K^{7+\ep}$ matrices in $\M^0_3(\K)$ whose lower right-hand entry $m_{33}$ is nonzero.

If any of the entries in the last row of $M$ is nonzero, then we can permute the columns of $M$ to bring that entry into the lower right-hand position; each such resulting matrix corresponds to at most three matrices in $\M^0_3(\K)$, and so there are still $\ll_\ep \K^{7+\ep}$ matrices in $\M^0_3(\K)$ that have any nonzero entry in the last row. Finally, any matrix whose last row consists of all zeros is certainly in $\M^0_3(\K)$, but there are only $(2\K+1)^6 \ll \K^6$ such matrices. We conclude that the total number of matrices in $\M^0_3(\K)$ is $\ll_\ep \K^{7+\ep}$, so that the probability of a randomly chosen matrix from $\M_3(\K)$ being singular is $|\M^0_3(\K)|/|\M_3(\K)| \ll_\ep  \K^{-2+\ep}$ as claimed.

{\bf Inductive step for $n\ge4$}:
Write a matrix $M\in\M_n(\K)$ in the form \eqref{Z.form}. Some such matrices will have $\det Z=0$; however, by the induction hypothesis for $n-2$, the probability that this occurs is $\ll_{n,\ep} \K^{-2+\ep}$ (independent of the entries outside $Z$), which is an allowably small probability.

Otherwise, fix values in $\I_\K$ for the $n^2-4$ entries other than $m_{11}$, $m_{12}$, $m_{21}$, and $m_{22}$ such that $\det Z\ne 0$.  It suffices to show that conditioning on any such fixed values, the probability that $M$ is singular, as $m_{11}$, $m_{12}$, $m_{21}$, and $m_{22}$ range over $\I_\K$, is $\ll_{n,\ep} \K^{-2+\ep}$.

By Lemma \ref{determinant.identity}, we see that $\det M=0$ is equivalent to $a_{11}a_{22} = a_{12}a_{21}$.  Recall that $a_{11}$ is a linear polynomial in the variable $m_{22}$ with leading coefficient $\det Z$, while the cofactor $a_{22}$ is a linear polynomial in $m_{11}$ with leading coefficient $\det Z$ (and similarly for the pair $a_{12}$ and $a_{21}$).  Moreover, the coefficients of these linear forms are sums of at most $(n-1)!$ products of $n-1$ entries at a time from $M$, hence are $\ll_n \K^{n-1}$ in size.  We may thus apply Lemma~\ref{not.many.solutions} with $\alpha=n-1$ to see that the probability of $a_{11}a_{22} = a_{12}a_{21}$ is $\ll_{n,\ep} \K^{-2+\ep}$, as desired. \qed
\medskip

Having established a suitably strong upper bound for $|\M^0_n(\K)|$, we can also bound the cardinality of $\M^\lambda_n(\K)$ for a fixed $\lambda \in \Z$,
using the fact that $M-\lambda I_n$ will be a singular matrix.
Notice that $\lambda$ can be as large as $n\K$ if we take $M$ to be the $n\times n$ matrix with all entries equal to $\K$, or as small as $-n\K$ if we take $M$ to be the $n\times n$ matrix with all entries equal to $-\K$.  It is not hard to show that these are the extreme cases for integer eigenvalues, and in fact even more is true:

\begin{lemma}
If $M \in \M_n(\K)$, then every complex eigenvalue of $M$ is at most $n\K$ in modulus.
\label{baby.Gershgorin}
\end{lemma}

\begin{proof}
Let $\lambda$ be any eigenvalue of $M$, and let $\v$ be a corresponding eigenvector, scaled so that $\max_{1\le i\le n} |\v_i| = 1$ (this is possible since $\v\ne\mathbf{0}$). Then
\[
|\lambda| = \max_{1\le i\le n} |(\lambda\v)_i| = \max_{1\le i\le n} |(M\v)_i| = \max_{1\le i\le n} \bigg| \sum_{k=1}^n m_{ik}\v_k \bigg|.
\]
Since each entry of $M$ is at most $\K$ in absolute value, and each coordinate of $\v$ is at most 1 in absolute value, we deduce that
\[
|\lambda| \le \max_{1\le i\le n} \sum_{k=1}^n |m_{ik}\v_k| \le \max_{1\le i\le n} \sum_{k=1}^n \K = n\K.
\]
\end{proof}

\begin{remark}
If we use the notation $D(z,r)$ to denote the disk of radius $r$ around the complex number $z$, then Lemma \ref{baby.Gershgorin} is the statement that every eigenvalue of $M$ must lie in $D(0,n\K)$. We remark that this statement is a weaker form of Gershgorin's ``circle theorem'' \cite{G}, which says that all of the eigenvalues of $M$ must lie in the union of the disks
\[
D\bigg( m_{11}, \sum_{\substack{1\le j\le n \\ j\ne1}} |m_{1j}| \bigg), \quad D\bigg( m_{22}, \sum_{\substack{1\le j\le n \\ j\ne2}} |m_{2j}| \bigg), \quad \dots, \quad D\bigg( m_{nn}, \sum_{\substack{1\le j\le n \\ j\ne n}} |m_{nj}| \bigg).
\]
In fact the proof of Gershgorin's theorem is very similar to the proof of Lemma \ref{baby.Gershgorin}, except that one begins with $|\lambda-m_{ii}|$ on the left-hand side rather than $|\lambda|$, and the other entries of $M$ are left explicit in the final inequality rather than estimated by $\K$ as above.
\end{remark}

We now have everything we need to prove the Theorem stated earlier: {\em given any integer $n\ge2$ and any real number $\ep>0$, the probability that a randomly chosen matrix in $\M_n(\K)$ has an integer eigenvalue is $\ll_{n,\ep} 1/\K^{1-\ep}$.}
By Lemma \ref{baby.Gershgorin}, any such integer eigenvalue $\lambda$ is at most $n\K$ in absolute value.
For each individual $\lambda$, we observe that if $M\in\M_n(\K)$ has eigenvalue $\lambda$, then $M-\lambda I_n$ is a singular matrix with integer entries which are bounded in absolute value by $\K+|\lambda| \le (n+1)\K$. Therefore every matrix in $\M^{\lambda}_n(\K)$ is contained in the set
\[
\big\{ M+\lambda I_n\colon M \in \M^0_n\big((n+1)\K\big)\big\}.
\]
By Lemma~\ref{any.dimension}, the cardinality of this set is $\ll_{n,\ep} ((n+1)\K)^{n^2-2+\ep} \ll_{n,\ep} \K^{n^2-2+\ep}$ for any fixed $\lambda$.  Summing over all values of $\lambda$ between $-n\K$ and $n\K$ (admittedly, some matrices are counted multiple times, but the upper bound thus obtained is still valid), we conclude that the total number of matrices in $\M^\Z_n(\K)$ is $\ll_{n,\ep} \K^{n^2-1+\ep}$. In other words, the probability that a matrix in $\M_n(\K)$ has an integer eigenvalue is $|\M^\Z_n(\K)|/|\M_n(\K)| \ll_{n,\ep} 1/\K^{1-\ep}$, as desired.

\section{More exact results for $2\times2$ matrices}
\label{2by2}

In the special case $n=2$, we can sharpen Lemma~\ref{any.dimension} and the Theorem considerably.  The $2\times 2$ case is particularly nice: since the trace of a matrix with integer entries is itself an integer, it follows that if one eigenvalue is an integer then both are. Consequently there is no distinction between being diagonalizable over $\Q$ and belonging to $\M^\Z_2(\K)$.

Although establishing these sharper results uses only standard techniques from analytic number theory, the computations are lengthy and require some tedious case-by-case considerations. Therefore we will content ourselves with simply giving the formulas in this section (proofs will appear in a subsequent paper~\cite{MW}). Using the notation $f(k) \sim g(k)$, which means that $\lim_{k\to\infty} f(k)/g(k) = 1$, we can state the two formulas in the following way:
The probability that a matrix in $\M_2(\K)$ is singular is asymptotically
\begin{equation}
\frac{|\M^0_2(\K)|}{|\M_2(\K)|} \sim \frac{6}{\pi^2} \left(\frac{\log \K}{\K^2}\right),
\label{sharper lemma}
\end{equation}
while the probability that a matrix in $\M_2(\K)$ has an integer eigenvalue is asymptotically
\begin{equation}
\frac{|\M^{\Z}_2(\K)|}{|\M_2(\K)|} \sim \left(\frac{7\sqrt{2} + 4 + 3\log(\sqrt{2}+1)}{3 \pi^2}\right)\left(\frac{\log \K}{\K}\right).
\label{sharper theorem}
\end{equation}
The orders of magnitude $(\log k)/k^2$ and $(\log k)/k$ are sharpenings of the orders of magnitude $k^{-2+\ep}$ in Lemma~\ref{any.dimension} and $k^{-1+\ep}$ in the Theorem, respectively.

A consequence of the asymptotic formulas~\eqref{sharper lemma} and~\eqref{sharper theorem} is that if we choose, uniformly at random, only matrices from $\M_2(\K)$ {\em with integer eigenvalues}, then the probability that such a matrix is singular is asymptotic to $4\alpha/\K$, where we have defined the constant
\begin{equation}
\alpha = \frac9{14\sqrt{2}+8+6\log(\sqrt{2}+1)} = 0.272008\dots.
\label{definition of alpha}
\end{equation}
In other words, the normalized quantity $\K |\M^0_2(\K) | / |\M^\Z_2(\K)|$ converges to $4\alpha$ as $\K$  tends to infinity. This convergence turns out to be a special case of a more general phenomenon:
if we hold $\lambda \in \Z$ fixed and let $\K$ tend to infinity, each normalized quantity $\K |\M^\lambda_2(\K)| / |\M^\Z_2(\K)|$ converges to the same constant $4\alpha$.

However, an interesting picture emerges if instead we rescale $\lambda$ along with $\K$, by thinking of $\lambda$ as the nearest integer to $\delta \K$ with $\delta$ a fixed real number. In fact, there is a continuous function $U^\Z(\delta)$ such that $\K |\M^\lambda_2(\K)| / |\M^\Z_2(\K)| \sim U^\Z(\lambda/\K)$ as $\K$ and $\lambda$ tend to infinity proportionally to each other.  The graph of $U^\Z(\delta)$ is the solid line in Figure~\ref{figure} below; the exact function is given by $U^\Z(\delta) = \alpha V(|\delta|)$, where $\alpha$ is the constant defined in equation~\eqref{definition of alpha} and
\begin{equation}
V(\delta) = 
\begin{cases}
4 - 2\delta - \delta^2 + \delta^2 \log(1+\delta) + 2(\delta-1)\log |\delta-1|,
& \text{if $0 \le \delta \le \sqrt{2}$},\\
\delta^2 - 2\delta - \log(\delta-1) - (\delta-1)^2 \log(\delta-1),
& \text{if $\sqrt{2} \le \delta \le 2$}.
\end{cases}
\end{equation}
(Note that there is no need to consider $V(\delta)$ for values of $|\delta|$ greater than 2, since all eigenvalues of matrices in $\M_2(\K)$ are at most $2k$ in modulus by Lemma~\ref{baby.Gershgorin}.)

\begin{figure}
\centering
\includegraphics[width=4in]{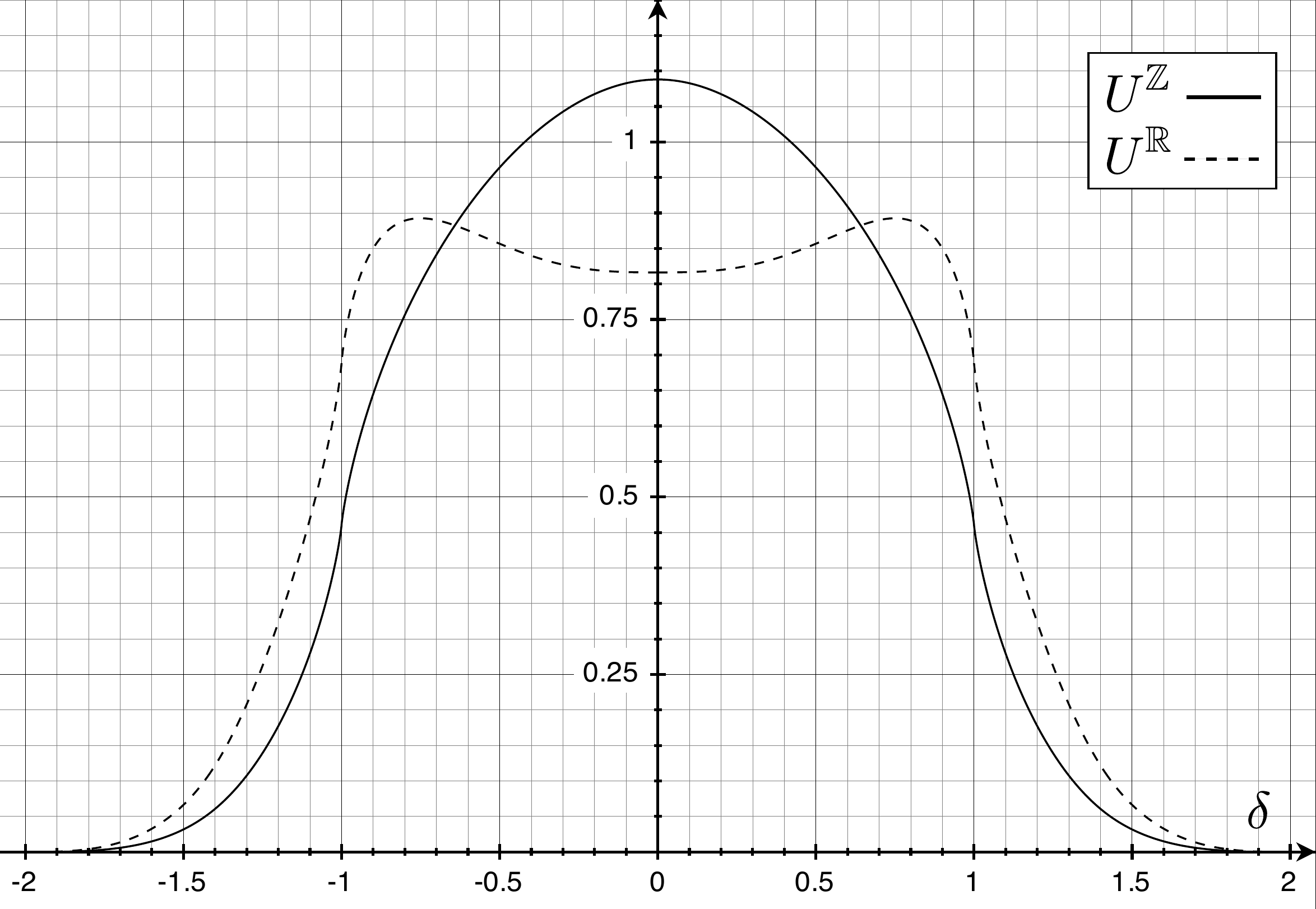}
\caption{Limiting distributions of real and integer eigenvalues for $\M_2(\K)$}
\label{figure}
\end{figure}

Intuitively, we can think the graph of $U^\Z(\delta)$ as follows.  For each positive integer $\K$, consider the histogram of eigenvalues from $\M^\Z_2(\K)$, vertically normalized by a factor of $1/|\M^\Z_2(\K)|$.  The total area under this histogram is exactly 2, since each matrix contributes exactly two eigenvalues.  If we then scale the horizontal axis by a factor of $1/\K$ and the vertical axis by a corresponding $\K$, the total area remains equal to~2, while the horizontal extent of the histogram lies in the interval $[-2,2]$. There is one such rescaled histogram for every positive integer $\K$; as $\K$ tends to infinity, the rescaled histograms converge pointwise to the limiting curve~$U^\Z(\delta)$.

(The astute reader will notice that we have ignored the fact that matrices with repeated eigenvalues occur only once in $\M^\lambda_2(\K)$ but contribute two eigenvalues to the histogram. This effect turns out to be negligible:  one can show by an argument similar to Lemma~\ref{not.many.solutions} that the number of matrices in $\M^\lambda_2(\K)$ with repeated eigenvalue $\lambda$ is $\ll_\ep \K^{1+\ep}$, a vanishingly small fraction of $|\M^\lambda_2(\K)|$).

For comparison, we can perform the exact same limiting process with the much larger subset of  $\M_2(\K)$ of matrices having real eigenvalues, which we naturally denote $\M^\R_2(\K)$.  In \cite{HLM}, Hetzel {\em et al} showed that the probability that a matrix in $\M_2(\K)$ has real eigenvalues, namely $|\M^\R_2(\K)|/|\M_2(\K)|$, converges to $49/72$ as $\K$ tends to infinity.  They observed that this probability can be realized as a Riemann sum for the indicator function of the set
\[
\{ (a,b,c,d)\in\R^4\colon |a|,|b|,|c|,|d|\le1,\, (a-d)^2+4bc \ge 0\},
\]
where the last inequality is precisely the condition for the matrix
$\big( \genfrac{}{}{0pt}{1}{a}{c} \genfrac{}{}{0pt}{1}{b}{d} \big)$
to have real eigenvalues.

If we likewise plot the histogram of eigenvalues from $\M^\R_2(\K)$, normalized to have area 2 as before, and scale horizontally by $1/\K$ and vertically by $\K$, we again get a limiting curve $U^\R(\delta)$ (which bounds an area of exactly~2 as well).  The graph of the function $U^\R(\delta)$ is the dashed line in Figure~\ref{figure}.  We can compute this curve using an integral representation similar to the one used to derive the constant $49/72$; although this integral is significantly more unwieldy, it eventually yields the exact formula $U^\R(\delta) = \beta W(|\delta|)$, where $\beta = 72/49$ and
\[
W(\delta) = 
\begin{cases}
(80 + 20\delta + 90\delta^2 + 52\delta^3 - 107\delta^4)/(144(1+\delta))
\\ \qquad{} - (5-7\delta+8\delta^2)(1-\delta) \log(1-\delta)/12
\\ \qquad{} - \delta(1-\delta^2) \log(1+\delta)/4,
  & \text{if $0 \le \delta \le 1$},\\
\delta(20+10\delta-12\delta^2-3\delta^3)/(16(1+\delta))
\\ \qquad{} + (3\delta-1)(\delta-1)\log(\delta-1)/4
\\ \qquad{} + \delta (\delta^2-1) \log(\delta+1)/4,
  & \text{if $1 \le \delta \le \sqrt{2}$},\\
\delta(\delta-2)(2-6\delta+3\delta^2)/(16(1+\delta))
\\ \qquad{}  - (\delta-1)^3 \log(\delta-1)/4,
  & \text{if $\sqrt{2} \le \delta \le 2$}.
\end{cases}
\]


It is interesting to note the qualitative differences between $U^\Z$ and $U^\R$.  Both are even functions, since $\M_2(\K)$ is closed under negation, and both functions are differentiable, even at $\delta = \pm \sqrt{2}$ and $\delta=\pm2$ (except for the points of infinite slope at $\delta =\pm 1$).  But the real-eigenvalue distribution is bimodal with its maxima at $\delta \approx \pm 0.7503$, while the integer-eigenvalue distribution is unimodal with its maximum at $\delta = 0$.  So a random $2\times 2$ matrix with integer entries bounded in absolute value by $\K$ is more likely to have an eigenvalue near $3\K/4$ than an eigenvalue near~0; but if we condition on having integer eigenvalues, then 0 becomes more likely than any other eigenvalue.

\smallskip
{\it Acknowledgments.} The authors were supported in part by grants from the Natural Sciences and Engineering Research Council.

\end{document}